\documentclass[12pt]{article}
\usepackage{amsmath,amsthm,amsfonts,latexsym,amssymb,fancyhea,supertab}
\usepackage{a4}

  \def\Q{\mathbb{Q}}
\def\Z{\mathbb{Z}} \def\R{\mathbb{R}}

\def\cE{{\cal E}}  \def\cF{{\cal F}}      

\def\al{\alpha}

\def\vn{\vspace{1ex}\noindent}
\def\vvn{\vspace{2ex}\noindent}

\newtheorem{theorem1}{Theorem}
\newtheorem{theorem2}[theorem1]{Theorem}
\title{Lucas sequences whose 12th or 9th term is a square}

\author{A.~Bremner\thanks{Department of Mathematics, Arizona State University, Tempe AZ, USA, 
e-mail: bremner@asu.edu , http://math.la.asu.edu/\~{}andrew/bremner.html} \and 
N.~Tzanakis\thanks{Department of Mathematics, University of Crete,
Iraklion, Greece, e-mail: tzanakis@math.uoc.gr , http://www.math.uoc.gr/\~{}tzanakis}
}

\date{}
\begin{document}

\maketitle
\section{Introduction}
      \label{introduction}
Let $P$ and $Q$ be non-zero relatively prime integers. The Lucas sequence 
$\{U_n(P,Q)\}$ is defined by
\begin{equation}
\label{Lucas}
U_0=0, \quad U_1=1, \quad U_n= P U_{n-1}-Q U_{n-2} \quad (n \geq 2). 
\end{equation}
The sequence $\{U_n(1,-1)\}$ is the familiar Fibonacci sequence, and it was proved by Cohn~
\cite{Co1} in 1964 that the only perfect square
greater than $1$ in this sequence is $U_{12}=144$. The question arises, for which parameters $P$, 
$Q$, can $U_n(P,Q)$ be a perfect square?
This has been studied by several authors: see for example Cohn~\cite{Co2}~\cite{Co3}~\cite{Co4}, 
Ljunggren~\cite{Lj}, and Robbins~\cite{Rob}. 
Using Baker's method on linear forms in logarithms, work of Shorey \& Tijdeman~\cite{ST} implies 
that there can only be finitely many squares in the sequence $\{U_n(P,Q)\}$. Ribenboim and 
McDaniel~\cite{RM1} with only elementary methods show that when $P$ and $Q$ are {\it odd}, and 
$P^2-4Q>0$, 
then $U_n$ can be square only for $n=0,1,2,3,6$ or $12$; and that there are at most two indices 
greater than 1 for which $U_n$ can be square. 
They characterize fully the instances when $U_n=\Box$, for $n=2,3,6$, and observe that 
$U_{12}=\Box$ if and only if there is a solution to the Diophantine system
\begin{equation}
\label{Paulocurve}
P=\Box, P^2-Q = 2\Box, P^2-2Q=3\Box, P^2-3Q=\Box, (P^2-2Q)^2-3Q^2=6\Box.
\end{equation}
When $P$ is {even}, a later paper of Ribenboim and McDaniel~\cite{RM2} proves that if 
$Q \equiv 1 \pmod 4$, then $U_n(P,Q)=\Box$ for $n>0$
only if $n$ is a square or twice a square, and all prime factors of $n$ divide $P^2-4Q$. 
Further, if $p^{2t}|n$ for a prime $p$, then
$U_{p^{2u}}$ is square for $u=1,\ldots,t$. In addition, if $n$ is even, then $U_n = \Box$ only if 
$P$ is a square or twice a square. 
A remark is made that no example is known of an integer pair $P$, $Q$, and an odd prime $p$, 
such that $U_{p^2}=\Box$ (note none can exist for $P$, $Q$ odd, $P^2-4Q>0$).

\vvn
In this paper, we complete the results of Ribenboim and MacDaniel~\cite{RM1} by determining 
all Lucas sequences $\{U_n(P,Q)\}$ with $U_{12}=\Box$ (in fact, the result is extended, because 
we do not need the restrictions that $P$, $Q$ be odd, and $P^2-4Q>0$): it turns 
out that the Fibonacci sequence provides the only example. Moreover, we also determine all Lucas 
sequences $\{U_n(P,Q)\}$ with $U_9=\Box$, subject only to the restriction that $(P,Q)=1$. 
Throughout this paper the symbol $\Box$ means square of a {\em non-zero} rational number. 
\begin{theorem1} \label{Th1}
Let $(P,Q)$ be any pair of relatively prime non-zero integers.
Then,   
\begin{itemize}
\item $U_{12}(P,Q)=\Box$ iff $(P,Q)=(1,-1)$ (corresponding to the Fibonacci sequence). 
\item $U_9(P,Q)=\Box$ iff $(P,Q)=(\pm 2,1)$ (corresponding to the sequences $U_n=n$ and 
$U_n=(-1)^{n+1}n$).
\end{itemize}
\end{theorem1}
The remainder of the paper is devoted mainly to the proof of this theorem.
Theorems 3 and 6 of \cite{RM1} combined with the first statement of Theorem \ref{Th1} imply the 
following.  
\begin{theorem2} \label{Th2}
 Let $P,Q$ be relatively prime odd integers, such that $P^2-4Q>0$. Then the $n$th term, $n>1$, of 
the Lucas sequence 
$U_n=U_n(P,Q)$ can be a square only if $n=2,3,6$ or $12$. More precisely\footnote{Below it is 
understood that parameters $a,b$ are 
in every case chosen so that $P,Q$ are odd, relatively prime and $P^2-4Q>0$.}:
\begin{itemize}
\item $U_2=\Box$ iff $P=a^2$.
\item $U_3=\Box$ iff $P=a$, $Q=a^2-b^2$.
\item $U_6=\Box$ iff $\displaystyle{P=3a^2b^2\,,\,Q=\frac{-a^8+12a^4b^4-9b^8}{2}}$.
\item $U_{12}=\Box$ iff $(P,Q)=(1,-1)$. Moreover, this result is also valid even  if we remove 
all restrictions on $P,Q$ except for $\gcd(P,Q)=1$.
\end{itemize}
\end{theorem2}

The proof of Theorem \ref{Th1} hinges, in both cases, upon finding all rational points on a curve 
of genus $2$. When the rank of the Jacobian of such a curve
is less than $2$, then methods of Chabauty~\cite{Ch}, as expounded subsequently by 
Coleman~\cite{Col}, Cassels and Flynn \cite{Ca-Fl} and 
Flynn \cite{F1} may be used to determine the (finitely many) rational points on the curve. When 
the rank of the Jacobian is at least $2$, as is the case here, a direct application of these 
methods fails. In order to deal with such situations, very interesting methods have been 
developed recently by a number of authors;  
see chapter 1 of Wetherell's Ph.D. thesis~\cite{W}, Bruin~\cite{Br1,Br2,Br3}, 
Bruin and Flynn~\cite{Br-Fl1,Br-Fl2}, Flynn~\cite{F2}, and Flynn \& Wetherell~\cite{FW1,FW2}.
For the purpose of this paper, the method of  \cite {F2} or \cite{FW1} is sufficient. 
\section{The Diophantine equations}
\label{DioEqns}
\subsection{The case $U_{12}$}
\label{u12genus2}
For $U_{12}(P,Q)$ to be square, we have from (\ref{Lucas})
\begin{equation}
\label{u12eqn}
U_{12}(P,Q)=P (P^2-3Q)(P^2-2Q)(P^2-Q)(P^4-4 P^2 Q+Q^2) = \Box.
\end{equation}
Now $(P(P^2-3Q)(P^2-Q), (P^2-Q)(P^4-4 P^2 Q+Q^2))$ divides $2$,
so that $U_{12}=\Box$ implies
\[ P(P^2-3Q)(P^2-Q) = \delta \Box, \quad (P^2-2Q)(P^4-4 P^2 Q+Q^2) = \delta \Box, \]
where $\delta = \pm 1, \pm 2$.
With  $x = Q/P^2$, we deduce
\[(1-2 x)(1-4 x+x^2) = \delta \Box, \]
and of these four elliptic curves, only the curve with $\delta=2$
has positive rational rank. Torsion points on the three other curves
do not provide any solutions for $P$, $Q$.
We are thus reduced to considering the equations
\[ P(P^2-3Q)(P^2-Q) = 2 \Box, \quad (P^2-2Q)(P^4-4 P^2 Q+Q^2) = 2 \Box. \]
From the first equation,
\[ P(P^2-3Q) = \pm 2 \Box, P^2-Q = \pm \Box, \quad \mbox{or} \quad  P(P^2-3Q) = \pm \Box,
                       P^2-Q = \pm 2 \Box. \]
The former case implies one of
\begin{eqnarray}
\label{array1}
P=\delta \Box & P^2-3Q=2\delta \Box & P^2-Q=\Box \nonumber \\
P=\delta \Box & P^2-3Q=-2\delta \Box & P^2-Q=-\Box \nonumber \\
P=2\delta \Box & P^2-3Q=\delta \Box & P^2-Q=\Box \nonumber \\
P=2\delta \Box & P^2-3Q=-\delta \Box & P^2-Q=-\Box
\end{eqnarray}
where $\delta=\pm 1, \pm 3$. \\
The latter case implies one of
\begin{eqnarray}
\label{array2}
P=\delta \Box & P^2-3Q=\delta \Box & P^2-Q=2\Box \nonumber \\
P=\delta \Box & P^2-3Q=-\delta \Box & P^2-Q=-2\Box
\end{eqnarray}
where $\delta=\pm 1, \pm 3$. \\
Solvability in $\R$ or elementary
congruences shows impossibility of the above equations (\ref{array1}), (\ref{array2}), except in 
the following instances:
\begin{eqnarray}
\label{array3}
P=-\Box, &  P^2-3Q=-2\Box, & P^2-Q=\Box \nonumber \\
P=-3\Box, & P^2-3Q=-6\Box, & P^2-Q=\Box \nonumber \\
P=6\Box, & P^2-3Q=3\Box, & P^2-Q=\Box \nonumber \\
P=\Box, & P^2-3Q=-2\Box, & P^2-Q=-\Box \nonumber \\
P=6\Box, & P^2-3Q=-3\Box, & P^2-Q=-\Box \nonumber \\
P=\Box, & P^2-3Q=\Box, & P^2-Q=2\Box \nonumber \\
P=-3\Box, & P^2-3Q=-3\Box, & P^2-Q=2\Box.
\end{eqnarray}
Recall now that
\[ (P^2-2Q)(P^4-4 P^2 Q+Q^2) = 2 \Box, \]
from which
\[ P^2-2Q = \eta \Box, P^4-4 P^2 Q+Q^2 = 2\eta \Box \quad \mbox{or} \quad P^2-2Q 
                                        = 2\eta \Box, P^4-4 P^2 Q+Q^2 = \eta \Box, \]
where $\eta=\pm 1, \pm 3$.
The only locally solvable equations are
\begin{eqnarray}
\label{array4}
P^2-2Q = -\Box, & P^4-4 P^2 Q+Q^2 = -2\Box \nonumber \\
P^2-2Q = 3\Box, & P^4-4 P^2 Q+Q^2 = 6\Box \nonumber \\
P^2-2Q = 2\Box, & P^4-4 P^2 Q+Q^2 = \Box \nonumber \\
\end{eqnarray}
It is straightforward by elementary congruences to deduce
from (\ref{array3}), (\ref{array4}), that we must have one of the following:
\begin{center}
\begin{tabular}{|r|r|r|r|r|}  \hline
$P$ & $P^2-3Q$ & $P^2-Q$ & $P^2-2Q$ & $P^4-4P^2Q+Q^2$ \\ \hline
$-\Box$ & $-2\Box$ & $\Box$ & $-\Box$ & $-2\Box$  \\ \hline
$6\Box$ & $3\Box$ & $\Box$ & $2\Box$ & $\Box$ \\ \hline
$\Box$ & $-2\Box$ & $-\Box$ & $-\Box$ & $-2\Box$ \\ \hline
$\Box$ & $\Box$ & $2\Box$ & $3\Box$ & $6\Box$ \\ \hline
\end{tabular}
\end{center}
Now the rational ranks of the following elliptic curves are 0:
\[ (-x+1)(x^2-4x+1)=-2\Box, \quad (-3x+1)(x^2-4x+1)=3\Box, \quad (-x+1)(x^2-4x+1)=2\Box, \]
and consequently the rational points on the curves 
corresponding to the first three rows of the above table
are straightforward to determine: they are $(P,Q)=(-1,1)$, $(0,-1)$, and
$(1,1)$ respectively. These lead to degenerate Lucas sequences with $U_{12}=0$.\\
It remains only to find all rational points on the following curve:
\[P=\Box, P^2-3Q=\Box, P^2-Q=2\Box, P^2-2Q=3\Box, P^4-4 P^2 Q+Q^2 =6\Box\,, \]
satisfying $(P,Q)=1$.
Note that this is the curve (\ref{Paulocurve}), though we have 
removed the restriction that
$P$ and $Q$ be odd, and $P^2-4Q>0$.\\
Put $Q/P^2 = 1-2 u^2$, so that 
\begin{equation}
\label{genus9}
3u^2-1 = 2\Box, \qquad 4u^2-1 = 3\Box, \qquad 2u^4+2u^2-1 = 3\Box.
\end{equation}
The equations (\ref{genus9}) define a curve of genus $9$, with
certainly only finitely many points.
We restrict attention to the curve of genus $2$ defined by
\[ 4u^2-1 = 3\Box, \qquad 2u^4+2u^2-1 = 3\Box. \]
Define $K=\Q(\sqrt{3})$, with ring of integers $\mathcal{O}_K=\Z[\sqrt{3}]$,
and fundamental unit $2+\sqrt{3}$.
Observe that $(u^2-1)^2-3 u^4 = - 3\Box$ implies
\[u^2 -1 + u^2 \sqrt{3} = \epsilon \sqrt{3} \gamma^2, \]
for $\epsilon$ a unit of $\mathcal{O}_K$ of norm $+1$,
and $\gamma \in \mathcal{O}_K$.
If $\epsilon=2+\sqrt{3}$, the resulting
equation is locally unsolvable above $3$, and so without loss of 
generality, $\epsilon=1$. Consider now
\[ u^2 (4u^2-1) (u^2(1+\sqrt{3})-1 ) =  3 \sqrt{3} V^2, \quad V \in K. \]
In consequence, $(x,y)=((12+4\sqrt{3})u^2, (36+12\sqrt{3})V)$ 
is a point defined over $K$ on the elliptic curve
\begin{equation}
\label{E_1}
E_1:  y^2 = x (x-(3+\sqrt{3})) (x-4\sqrt{3})
\end{equation}
satisfying $\frac{(3-\sqrt{3})}{24} x \in \Q^2$.
We shall see that the $K$-rank of $E_1$ is equal to $1$, 
with generator of infinite order $P=(\sqrt{3},3\sqrt{3})$.
\subsection{The case $U_9$}
For $U_9(P,Q)$ to be square, we have from (\ref{Lucas})
\begin{equation}
\label{u9eqn}
U_9(P,Q) = (P^2-Q)(P^6-6 P^4 Q+9 P^2 Q^2-Q^3) = \Box.
\end{equation}
So
\[ P^2-Q = \delta \Box, \qquad P^6-6 P^4 Q+9 P^2 Q^2-Q^3 = \delta \Box, \]
where $\pm \delta = 1, 3$. Put $Q=P^2-\delta R^2$. Then
\begin{equation}
\label{u9eq}
\frac{3}{\delta} P^6-9 P^4 R^2+6 \delta P^2 R^4+\delta^2 R^6 = \Box, 
\end{equation}
with covering elliptic curves
\begin{equation}
\label{J1}
\frac{3}{\delta}-9 x +6 \delta x^2+\delta^2 x^3 = \Box, \qquad 
                  \frac{3}{\delta} x^3-9 x^2+6\delta x+\delta^2 = \Box.
\end{equation}
For $\delta=\pm 1$, the first curve has rational rank $0$, and torsion points
do not lead to non-zero solutions for $P$, $Q$.
For $\delta=\pm 3$, both curves at (\ref{J1}) have rational rank $1$,
so that the rank of the Jacobian of (\ref{u9eq}) equals $2$.
To solve the equation $U_9(P,Q)=\Box$, it is necessary to
determine all integer points on the two curves
\begin{equation}
\label{u91}
 P^6-9 P^4 R^2+18 P^2 R^4+9 R^6 = \Box,
\end{equation}
and
\begin{equation}
\label{u92}
-P^6-9 P^4 R^2 -18 P^2 R^4+9 R^6 = \Box. 
\end{equation}
To this end, let $L=\Q(\al)$ be the number field defined by 
$\al^3-3\al-1=0$. Gal($L/\Q$) is cyclic of order $3$, generated
by $\sigma$, say, where $\al^\sigma=-1/(1+\al)=-2-\al+\al^2$.
The ring of integers $\mathcal{O}_L$ has basis $\{1,\al,\al^2\}$,
and class number $1$. Generators for the group of units in 
$\mathcal{O}_L$ are $\epsilon_1 =\al, \epsilon_2=1+\al$, with norms
$\mbox{Norm}(\epsilon_1)=1$, $\mbox{Norm}(\epsilon_2)=-1$.
The discriminant of $L/Q$ is 81, and the ideal $(3)$ factors in $\mathcal{O}_L$
as $(-1+\al)^3$.
\subsubsection{}
\label{u9genus2}
Equation (\ref{u91}) may be written in the form
\[ \mbox{Norm}_{L/\Q} (P^2 +(-5+\al+\al^2) R^2) = S^2, \mbox{ say}, \]
and it follows that
\begin{equation}
\label{normeq1}
 P^2 +(-5+\al+\al^2) R^2 = \lambda U^2, 
\end{equation}
with $\lambda \in \mathcal{O}_L$ squarefree and of norm $+1$ modulo
$L^{*^2}$.\\
Applying $\sigma$,
\begin{equation}
\label{normeq2}
 P^2 +(-5-2\al+\al^2)R^2 = \lambda^\sigma V^2.
\end{equation}
Suppose $\mathcal{P}$ is a first degree prime ideal of
$\mathcal{O}_L$ dividing $(\lambda)$. Then for the norm of $\lambda$
to be a square, $\lambda$ must also be divisible by one of the conjugate
prime ideals of $\mathcal{P}$. It follows that $\mathcal{P}$, or one of
its conjugates, divides both $\lambda$ and $\lambda^\sigma$. Then
this prime will divide
$((-5+\al+\al^2)-(-5-2\al+\al^2)) = (3 \al) = (1-\al)^3$.
So $\mathcal{P}$ has to be $(1-\al)$, with $(1-\al)^2$ dividing
$\lambda$, contradicting $\lambda$ squarefree.
If the residual degree of $\mathcal{P}$ is 3, then the norm
of $\lambda$ cannot be square.
Finally, the residual degree of $\mathcal{P}$ cannot be 2, otherwise
$\theta = 5-\al-\al^2 \equiv m \pmod{\mathcal{P}}$, for
some rational integer $m$, so that $\al = 1-2\theta+\frac{1}{3}\theta^2$
is congruent to a rational integer modulo $\mathcal{P}$, impossible. 
In consequence, $\lambda$ is
forced to be a unit, of norm $+1$.  Without loss of generality, the only 
possibilities are
$\lambda=1, \epsilon_1, -\epsilon_2, -\epsilon_1 \epsilon_2.$
However, specializing the left hand side of (\ref{normeq1}) at the root
$\al_0=1.8793852415...$ of $x^3-3 x-1 = 0$ shows that
$P^2+ 0.4114..R^2 = \lambda(\al_0) U(\al_0)^2$, so that
$\lambda(\al_0) > 0$, giving unsolvability of (\ref{normeq1}) for
$\lambda = -\epsilon_2$, $-\epsilon_1 \epsilon_2$.
There remain the two cases $\lambda=1$, with solution $(P,R,U)=(1,0,1)$, and
$\lambda=\epsilon_1$, with solution $(P,R,U) = (0,1, 4-\al^2)$.
From (\ref{normeq1}) and (\ref{normeq2}) we now have
\[ P^2 (P^2+(-5+\al+\al^2)R^2) (P^2 +(-5-2\al+\al^2)R^2) = \mu W^2, \]
with $\mu = \lambda \lambda^\sigma = 1$ or $1+\al-\al^2$.
Accordingly, $X=P^2/R^2$ gives a point on the elliptic curve:
\begin{equation}
\label{E1}
 X (X+(-5+\al+\al^2)) (X+(-5-2\al+\al^2)) = \mu Y^2.
\end{equation}
Now when $\mu=1$, a relatively straightforward $2$-descent argument
shows that the $\Q(\al)$-rank of (\ref{E1}) is equal to $0$ (we also 
checked this result using the Pari-GP software of 
Denis Simon~\cite{Sim}). 
The torsion group is of order $4$, and no non-zero $P$, $Q$ arise.\\
When $\mu=1+\al-\al^2$, then 
$(x,y)=(\frac{\mu}{(1-\al)^2} \frac{P^2}{R^2}, \frac{\mu^2}{(1-\al)^3} \frac{W}{R^3})$
is a point on the elliptic curve $E_2$ over $\Q(\al)$, where
\begin{equation}
\label{E_2}
E_2 : y^2 = x(x+(-2-\al+\al^2)) (x+(-1+\al+\al^2)),
\end{equation}
satisfying
$\frac{(1-\al)^2}{\mu}x = (4+\al-2\al^2) x \in \Q^2$.

We shall see that the
$\Q(\al)$-rank of $E_2$ is $1$, with generator of infinite order
equal to
$(1,\al)$.
\subsubsection{}
Equation (\ref{u92}) may be written in the form
\[ \mbox{Norm}_{L/{\Q}} (-P^2 +(-5+\al+\al^2) R^2) = S^2, \mbox{ say}, \]
so that
\begin{equation}
\label{normeq3}
 -P^2 +(-5+\al+\al^2) R^2 = \lambda U^2,
\end{equation}
with $\lambda \in \mathcal{O}_L$ squarefree and of norm $+1$ modulo
$L^{*^2}$. Arguing as in the previous case, $\lambda$ must be a squarefree
unit of norm $+1$, so without loss of generality equal to $1$, $\epsilon_1$, $-\epsilon_2$, 
$-\epsilon_1 \epsilon_2$.  Only when $\lambda=\epsilon_1$ is (\ref{normeq3}) 
solvable at all the infinite places. Thus
\[ -P^2 +(-5+\al+\al^2) R^2 = \al U^2, \quad -P^2 +(-5-2\al+\al^2) R^2 = (-2-\al+\al^2) V^2, \]
and 
$x=\frac{(1+\al-\al^2)}{(1-\al)^2} \frac{P^2}{R^2}$
is the $x$-coordinate of a point on the elliptic curve 
\begin{equation}
\label{E2'}
y^2 = x (x + (2+\al-\al^2)) (x + (1-\al-\al^2)).
\end{equation}
satisfying 
$\frac{(1-\al)^2}{(1+\al-\al^2)}x=(4+\al-2\al^2)x \in \Q^2$.
However, a straightforward calculation shows that the $\Q(\al)$-rank
of $(\ref{E2'})$ is equal to $0$, with torsion group the obvious
group of order $4$. There are no corresponding solutions for $P$, $Q$.

\section{The Mordell-Weil basis} \label{MWbasis}
Here we justify our assertions about the elliptic curves $E_1$ at (\ref{E_1})
and $E_2$ at (\ref{E_2}). 
These curves are defined over fields $F$ (where $F$=$K$ or $L$, respectively)
with unique factorization, and have $F$-rational two-torsion. So standard two-descents
over $F$ work analogously to the standard two-descent over $\Q$ for an elliptic curve
with rational two-torsion; see for example Silverman~\cite{Sil}, Chapter 10.4. 
It is thus straightforward
to determine generators for $E_i(F)/2E_i(F)$, $i=1,2$ (and in fact software packages such
as that of Simon~\cite{Sim} written for Pari-GP, and ALGAE~\cite{BrA} of Bruin written for
KASH, with m-ALGAE~\cite{BrmA} for MAGMA, also perform this calculation effortlessly).
Such generators are the classes of $P_1=(\sqrt{3}, 3\sqrt{3})$ for the curve $E_1$, and 
$P_2=(1,\alpha)$
for the curve $E_2$. In fact $P_1$ and $P_2$ are generators for the Mordell-Weil groups
$E_1(\Q(\sqrt{3}))$ and $E_2(\Q(\alpha))$ respectively. To show this necessitates detailed height
calculations over the appropriate number field, with careful estimates for the difference
$\hat{h}(Q) - \frac{1}{2} h(Q)$ where $\hat{h}(Q)$ is the canonical height of the point $Q$,
and $h(Q)$ the logarithmic height. The KASH/TECC package of Kida~\cite{Ki} was
useful here in confirming calculations. The standard Silverman bounds~\cite{Sil2} are numerically
too crude for our purposes, so recourse was made to the refinements of Siksek~\cite{Sik}.
Full details of the argument are given in Section 3 of \cite{BT}.
Actually, determination of the full $F$-rational Mordell-Weil groups of $E_1$ and $E_2$
may be  redundant; it is likely that the subsequent
local computations can be performed subject only to a simple condition on
the index in $E(F)$ of a set of generators for $E(F)/2E(F)$.
The reader is referred
to Bruin~\cite{Br3} or Flynn \& Wetherell~\cite{FW2} for details and examples. This latter
technique must be used of course when the height computations are simply too time consuming
to be practical.

\section{Finding all points on (\ref{E_1}) and (\ref{E_2}) under the rationality conditions} 
\label{the solutions}
\subsection{General description of the method}
         \label{formalgroup}
The problems to which we were led in section~\ref{DioEqns} are of the following shape. 

{\bf Problem:} Let 
\begin{equation} \label{gen_xy_elliptic}
\cE\,:\, y^2+a_1xy+a_3y=x^3+a_2x^2+a_4x +a_6 \;,
\end{equation}
be an elliptic curve defined over $\Q(\alpha)$, where $\alpha$ is a root 
of a polynomial $f(X)\in\Z[X]$, irreducible over $\Q$, of degree $d\geq 2$ 
and $\beta\in\Q(\alpha)$ an algebraic integer. Find all points $(x,y)\in\cE(\Q(\alpha))$ 
for which $\beta x$ is the square of a {\em rational number}.

\vn
For the solution of this type of problem we adopt the technique described and applied in 
Flynn and Wetherell~\cite{FW1}\footnote{But see also the references at the end of 
section \ref{introduction}.}. Several problems of this type have already been solved with a 
similar technique (besides \cite{FW1}, see
also \cite{Br3},\cite{Br-Fl2},\cite{F2},\cite{FW2}); therefore we content ourselves with a rather 
rough description of the employed method and refer the interested reader to section 4 of 
\cite{BT} for a detailed exposition.

We assume the existence of a rational prime $p$ with the following properties: 

(i) $f(X)$ is irreducible in $\Q_p[X]$. 

(ii) The coefficients of (\ref{gen_xy_elliptic}) are in $\Z_p[\alpha]$.

(iii) Equation (\ref{gen_xy_elliptic}) is a minimal Weierstrass equation for $\cE/\Q_p(\alpha)$ 
at the unique discrete 
valuation $v$ defined on $\Q(\alpha)$ with $v(p)=1$.  

(iv) $\beta\in\Q_p(\alpha)$ is a $p$-adic unit. 

\vn
We work with both (\ref{gen_xy_elliptic}) and the associated model 
$w=z^3+a_1zw+a_2z^2w+a_3w^2+a_4zw^2+a_6w^3 $,
which are related by means of the birational transformation \\  
$ (x,y) \mapsto (z,w)=(-x/y,-1/y) \;,\; (z,w) \mapsto (x,y)=(z/w,-1/w) .$ 

\vn
We also need the {\em formal group law} which is defined by means of two $p$-adically convergent 
power series 
$\cF(z_1,z_2)\in\Z[\alpha][[z_1,z_2]]$ (``sum") and $\iota(z)\in\Z[\alpha][[z]]$ (``inverse"), 
satisfying certain properties (see \S\,2,chapter IV of~\cite{Sil}).
These series can be explicitly calculated up to any presision, and the operations 
$(z_1,z_2) \mapsto \cF(z_1,z_2)$, $z\mapsto \iota(z)$ make $p\Z_p[\alpha]$ a group $\hat{\cE}$ 
(or, more precisely, $\hat{\cE}/\Z_p[\alpha]$), which is the {\em formal group} associated to  
$\cE/\Q_p(\alpha)$.
 
There is a group isomorphism between $\hat{\cE}$ and the subgroup of ${\cal E}(\Q_p(\alpha))$ 
consisting of those points $Q$ whose reduction $\bmod{p}$ is the zero point of the reduced 
$\bmod{p}$ curve, defined by \\
$z \mapsto Q$, where $Q=(z/w(z)\,,\,-1/w(z))$ if $z\neq 0$ and $Q={\cal O}$  if $z=0$, \\
with $w(z)$ a $p$-adically convergent power series, that can be explicitly calculated up to any 
$p$-adic 
precision. The inverse map is given by $z({\cal O})=0$ and for $Q\in{\cal E}(\Q_p(\alpha))$ 
different from $\cal O$
whose reduction $\bmod{p}$ is zero, $z(Q)=-\frac{x(Q)}{y(Q)}$. 

\vn
The remarkable property relating the functions $z$ and $\cF$ is that, for any points $Q_1,Q_2$ as 
$Q$ above, \vspace{-2mm}
\begin{equation} \label{z_sum}
z(Q_1+Q_2) = \cF(z(Q_1),z(Q_2))\:.
\end{equation}
With respect to $\cal E$, a {\em logarithmic} function $\log$ is defined on $p\Z_p[\alpha]$ and 
an {\em exponential} function $\exp$ is defined on $p^r\Z_p[\alpha]$, where $r=1$ if $p>2$ and 
$r=2$ if $p=2$. These functions are mutually inverse and can be explicitly calculated as $p$-adic 
power series up to any precision. Moreover, if $r$ is as above and $z_1,z_2 \in p^r\Z_p[\alpha]$, 
then \vspace{-2mm}
\begin{equation} \label{z_linform}
\mbox{$\log\cF(z_1,z_2)=\log z_1 +\log z_2$ and $\exp(z_1+z_2)=\cF(\exp z_1,\exp z_2)$.}  
\end{equation}
Suppose now we know a point such that $z(Q)\in p^r\Z_p[\alpha]$ and assume further that, for a 
certain specifically known point  
$P\in\Q(\alpha)\cap \cE(\Z_p[\alpha])$, or for $P={\cal O}$, we want to find all $n\in\Z$ for 
which $\beta x(P+nQ)$ is a rational number (or, more particularly, a square of a rational). 
According to whether $P$ is a finite point or $P={\cal O}$, we express $\beta x(P+nQ)$ or 
$1/\beta x(nQ)$ first as an element of $\Z_p[\alpha][[z(nQ)]]$ and then, using properties 
(\ref{z_sum}) and (\ref{z_linform}), as a sum
$ \theta_0(n)+\theta_1(n)\alpha+\cdots +\theta_{d-1}(n)\alpha^{d-1}$,
where each series $\theta_i(n)$ is a $p$-adically convergent power series in $n$ with 
coefficients in $\Z_p$, which can be explicitly calculated up to any 
desired $p$-adic precision.  
In order that this sum be a rational number we must have $\theta_i(n)=0$ for $i=1,\ldots, d-1$. 
At this point we use Strassman's Theorem\footnote{Theorem 4.1, in \cite{Ca}.}, which restricts 
the number of $p$-adic integer solutions $n$. If the maximum number of solutions implied by this 
theorem is equal to the number of solutions that we {\em actually} know, then we have explicitly 
all solutions. Sometimes, 
as in some instances of  the sections below, even Strassman's Theorem is not necessary. 

In the following two sections we apply the above method to equations (\ref{E_1}) and (\ref{E_2}). 
At the suggestion of the referee, we give only a few computational details; most of our 
computational results, including the explicit form of the functions $\cF(z_1,z_2), w(z),\log$ and 
$\exp$ with the required precision, can be found in section 5.1 of \cite{BT}.  

\subsection{Equation (\ref{E_1})} \label{solve_u12}
For this section, let $\al=\sqrt{3}$. We write (\ref{E_1}) as 
$\cE: y^2 = x^3-(3+5\alpha)x^2+12(1+\alpha)x $ and  
according to the discussion in section \ref{u12genus2} we must find all points$(x,y)$ on this 
curve, such that $\beta x = u^2 \in\Q^2$, 
where $\beta=(3-\alpha)/24$.\\
We work $p$-adically with $p=7$. According to section~\ref{MWbasis}, any point on 
$\cE(\Q(\alpha))$ is of the form $n_1P_1+T$, where 
$P_1=(\alpha, 3\alpha)$ and  $T\in\{{\cal O},(0, 0),(4\alpha, 0),(3+\alpha,0) \}$.
For $Q=11P_1$ we have $z(Q)\in 7\Z_7[\alpha]$, and any point of $\cE(\Q(\alpha))$ can be written 
in the form
$n_1P_1+T=(11n+r)P_1+T=nQ+P$, 
with $P=rP_1+T\,,\: -5\leq r \leq 5$ and $T$ a torsion point as above. There are 44 possibilities 
for $P$, one of which is $P={\cal O}$.

(i) Consider first the case when $P$ is one out of the 43 possible finite points. As noted in 
section \ref{formalgroup}, we are led to a relation
$\beta x(P+nQ)=\theta_0(n)+\theta_1(n)\alpha$, where the $7$-adically convergent series 
$\theta_i(n), i=0,1$ depend on $P$.  
In 35 out of the 43 cases, it turns out that $\theta_1(0)\not\equiv 0\pmod{7}$, which, in 
particular, implies that  
$\beta x(P+nQ)=\theta_0(n)+\theta_1(n)\alpha \equiv \theta_0(0)+\theta_1(0)\alpha \pmod{7}$ 
cannot be rational. The only cases that are not excluded in this way occur when $P$ is one of the 
following
points: $\pm 4P_1+(0,0),\pm 3P_1+(0,0), \pm P_1+(0,0), (0,0),(3+\alpha,0)$.

We deal with these cases as follows: If $P=\pm 4P_1+(0,0)$, then $\theta_0(0)=5$, a quadratic 
non-residue of $7$; therefore, whatever $n$ may be, 
$\beta x(P+nQ)=\theta_0(n)+\theta_1(n)\alpha$ cannot be the square of a rational number. In a 
completely analogous manner we exclude $P=\pm 3P_1+(0,0)$, since, in this case, $\theta_0(0)=6$. 

\noindent 
Next, consider $P=\pm P_1+(0,0)=(12+4\alpha,\pm(36+12\alpha))$. With the plus sign we compute
$\theta_1(n)=7\cdot 94 n+7^2\cdot 40 n^2 +7^3\cdot 6 n^3 +\cdots $, and with the minus sign, 
$\theta_1(n)=7\cdot 249 n+7^2\cdot 40 n^2 +7^3 n^3 +\cdots $. In both cases, if $n\neq 0$, then 
dividing out by $7n$ we are led to an impossible relation $\bmod{7}$; hence $n=0$ and 
$x=x(P+nQ)=x(P)=12+4\alpha$, which gives 
$u^2=\beta (12+4\alpha)=1$ and $u=1$.  

\noindent
If $P=(0,0)$, then we compute 
$\theta_1(n)=7^3\cdot 6889 n^2 +7^4\cdot 1733 n^4 +7^7\cdot 2n^6+7^7\cdot 5n^8+\cdots $ and if 
$n\neq 0$ we divide out by $7^3n^2$ and we are led to an impossible relation $\bmod{7}$. Thus, 
$n=0$, which leads to $x=0$ and $u=0$.
Finally, if $P=(3+\alpha,0)$, then $\theta_1(n)=7^2\cdot 288 n^2 +7^4n^4+\cdots $, forcing again 
$n=0$. Thus,
$x=3+\alpha$ and $u^2=\beta (3+\alpha)=1/4$, hence $u=1/2$.

(ii) Assume now that $P={\cal O}$. Then we have $1/\beta x(nQ)=\theta_0(n)+\theta_1(n)\alpha$, 
where 
$\theta_1(n)=7^2\cdot 244 n^2+7^4\cdot 2n^4+\cdots $. Since we are interested in finite points 
$(x,y)=nQ$, $n$ must be non-zero. Dividing out 
$\theta_1(n)=0$ by $7^2n^2$ we obtain an impossible equality. 

{\bf Conclusion.} The only points on (\ref{E_1}) satisfying $\beta x=u^2\in\Q^2$ are those with 
$x=12+4\alpha, 3+\alpha,0$, corresponding to $u=1,\frac{1}{2},0$. Only the first leads to a 
solution of (\ref{genus9}) and this leads to the solution $(P,Q)=(1,-1)$ of (\ref{u12eqn}) with 
$U_{12}(1,-1)=12^2$.

\subsection{Equation (\ref{E_2})} \label{solve_u9}
We write (\ref{E_2}) as $ y^2 = x^3+(-3+2\alpha^2)x^2+(2-\alpha^2)x $. 
According to the discussion in section \ref{u9genus2}, it suffices to find all points $(x,y)$ 
on this curve  
such that $\beta x\in\Q^2$, where $\beta=4+\alpha-2\alpha^2$.
We work $p$-adically with $p=2$. According to section \ref{MWbasis}, any point on 
$\cE(\Q(\alpha))$ is of the 
form $n_1P_1+T$, where $P_1=(1, \alpha)$ and 
$T\in\{{\cal O}\,,\,(1-\alpha-\alpha^2, 0)\,,\,(0, 0)\,,\,(2+\alpha-\alpha^2,0) \}$.
For the point $Q=4P_1$ we have $z(Q)\in 4\Z_2[\alpha]$ and we write any point on 
$\cE(\Q(\alpha))$ in the form
$ n_1P_1+T=(4n+r)P_1+T=nQ+P$, with $P=rP_1+T\,,\: r\in\{-1,0,1,2\}$ and $T$ a torsion point as 
above.
Therefore, there are 16 possibilities for $P$, one of which is $P={\cal O}$.

Working as in section \ref{solve_u12} we check that the only solutions $(x,y)$ such that 
$\beta x\in\Q^2$ are $(x,y)=(0,0),2P_1+(0,0),-2P_1+(0,0)$. To give an idea of how we apply 
Strassman's Theorem, let us consider the instance when $P=2P_1+T$ with $T=(0,0)$.  We compute 
$\theta_1(n) = 2^6\cdot 7 n+2^6\cdot 3 n^2+0\cdot n^3+0\cdot n^4+0\cdot n^5+2^7n^6(\cdots)$. By 
Strassman's Theorem, $\theta_1(n)=0$ can have at most two solutions in $2$-adic integers $n$. On 
the other hand, a straightforward computation shows that $\beta x(P+0\cdot Q)=4=\beta x(P-Q)$, 
which implies, in particular, that $\theta_1(0)=0=\theta_1(-1)$. Hence, $n=0,-1$ are the only 
solutions obtained for the above specific value of $P$, leading to the points
$(x,y)=P+0\,Q=2P_1+(0,0)$ and $(x,y)=P+(-1)Q=-2P_1+(0,0)$, both having 
$x=\frac{4}{3}(1-\alpha^2)$ and $\beta x=4$.  

{\bf Conclusion.} The only points on (\ref{E_2}) satisfying $\beta x\in\Q^2$ are those with
$x=0$ (leading to $P=0$), and 
$x=\frac{4}{3}(1-\alpha^2)$ giving successively (in the notation of section~\ref{u9genus2}) 
$\frac{P^2}{R^2}=4$ and $(P,Q)=(\pm 2,1)$, corresponding to degenerate Lucas sequences.

\section*{Acknowledgement}
We are grateful to the referee for suggestions which have greatly improved the presentation 
of this paper.

\end{document}